\documentclass[10pt]{amsart} 
\usepackage{pstricks, pst-node, amssymb}
\usepackage{epsfig,color}
\usepackage{graphicx}
\usepackage{amsmath, amsfonts}
\usepackage{amsthm,amssymb,amsbsy, latexsym, stmaryrd, amscd, xy} 
\usepackage[mathscr]{eucal}
\usepackage{t1enc}
\usepackage{fancyheadings} 
\usepackage{epic} 
\usepackage{eepic}

\psset{unit=1.1pt, arrowsize=4pt, linewidth=.7pt}
\psset{linecolor=blue}
\newgray{grayish}{.90}
\newrgbcolor{embgreen}{0 .5 0}
\def\vblack(#1, #2)#3{\cnode*[linecolor=black](#1, #2){3}{#3}}
\def\vwhite(#1,#2)#3{\cnode[linecolor=black,fillcolor=white,fillstyle=solid](#1,#2){3}{#3}}
\countdef\x=23
\countdef\y=24
\countdef\z=25
\countdef\t=26

\def\tbox(#1,#2)#3{
\x=#1 \y=#2
\multiply\x by 12
\multiply\y by 12
\z=\x \t=\y
\advance\z by 12
\advance\t by 12
\psline(\x,\y)(\x,\t)(\z,\t)(\z,\y)(\x,\y)
\advance\x by 6
\advance\y by 6
\rput(\x,\y){{\bf #3}}}

\newtheorem{prop}{Proposition} 
\newtheorem{lemma}[prop]{Lemma}

\newtheorem{theorem}[prop]{Theorem}

\newtheorem{conj}[prop]{Conjecture}

\theoremstyle{definition}

\newtheorem{example}[prop]{Example}

\newcommand{\ca}{\mathcal{A}}
\newcommand{\cd}{\mathcal{D}}
\newcommand{\si}{\sigma}

\def\dd{\kern.4ex\mbox{\raise.7ex\hbox{{\rule{.45em}{.12ex}}}}\kern.4ex}

\def\newMAH#1{%
\expandafter\def\csname #1\endcsname{\mathop{\mbox{{\sc#1}}}\nolimits}%
}

\newMAH{maj}

\def\newexpMAH#1{%
\expandafter\def\csname exp#1\endcsname{\mathop{\mbox{{\footnotesize\sc{#1}}}}\nolimits}%
}

\newexpMAH{maj}

\def\biw#1,#2,{%
\begin{pmatrix}
#1\cr
#2\cr
\end{pmatrix}
}

\def\ch#1,#2,{{#1\choose #2}}

\newcommand{\cls}{\mathcal{S}}

\def\newop#1{\expandafter\def\csname#1\endcsname{\mathop{\rm#1}\nolimits}}

\def\newscop#1{%
\expandafter\def\csname #1\endcsname{\mathop{\mbox{{\sc#1}}}\nolimits}%
}

\def\newcapop#1{%
\expandafter\def\csname #1\endcsname{\mathop{\mbox{{\sc#1}}}\nolimits}%
}

\newscop{inv}
\newscop{mak}
\newscop{makl}

\def\emm#1,{{\em #1}}

\def\bskip#1,{\vspace*{#1\baselineskip}}

\newop{des}
\newop{wex}

\newscop{wb}
\newscop{wt}

\newscop{db}
\newscop{dt}

\newscop{wextopsum}
\newscop{wexbotsum}
\newscop{destopsum}
\newscop{desbotsum}

\newscop{wextopset}
\newscop{wexbotset}
\newscop{destopset}
\newscop{desbotset}

\newscop{rembr}

\newcapop{cab}

\newcommand\T{\ensuremath{\mathcal{T}}}

\catcode`\@=11
\def\section{\@startsection{section}{1}%
 \z@{.7\linespacing\@plus\linespacing}{.5\linespacing}%
 {\normalfont\bfseries\scshape\centering}}
\def\subsection{\@startsection{subsection}{2}%
  \z@{.5\linespacing\@plus\linespacing}{.5\linespacing}%
  {\normalfont\bfseries\scshape}}

\def\subsubsection{\@startsection{subsubsection}{3}%
  \z@{.5\linespacing\@plus.7\linespacing}{-.5em}%
  {\normalfont\itshape}}
\catcode`\@=12
\title{A conjecture of Stanley on alternating permutations} 

\author{Robin Chapman and Lauren K. Williams} 

\address{
  Department of Mathematics, University of Bristol,
   Bristol, BS8 1TW, UK}

\address{
  Department of Mathematics, Harvard University, Cambridge, MA 02138}

\email{lauren@math.harvard.edu} 

\date{\today}

\begin{document}

\begin{abstract}
We give two simple proofs of a 
conjecture of Richard Stanley concerning the equidistribution
of derangements and alternating permutations with the maximal number
of fixed points.  
\end{abstract}

\maketitle \thispagestyle{empty}

\section{Introduction}

We write $[n]=\{1,\dots , n\}$ and 
$S_n$ for the set of permutations of 
$[n]$.
A permutation is {\it alternating} if 
$a_1 > a_2 < a_3 > a_4 < \dots $.  Similarly, define $w$ to be
{\it reverse alternating} if $a_1 < a_2 > a_3 < a_4 > \dots $.

In \cite{Stanley}, Richard Stanley used the theory of symmetric functions 
to enumerate various classes of alternating
permutations $w$ of $\{1, 2, \dots n\}$.
One class that he considered were alternating permutations $w$ with
a specified number of fixed points.  

Write $d_k(n)$ (respectively, $d^*_k(n)$) for the number of 
alternating (respectively, reverse alternating) permutations in 
$S_n$ with $k$ fixed points.

As observed in \cite{Stanley}, it is not hard to see that 
\begin{align*}
\max \{ k: d_k(n) \neq 0 \} &= \lceil n/2 \rceil, n \geq 4 \\
\max \{ k: d^*_k(n) \neq 0 \} &= \lceil (n+1)/2 \rceil, n \geq 5. 
\end{align*}

Stanley conjectured the following in \cite{Stanley}.

\begin{conj}\cite[Conjecture 6.3]{Stanley}\label{main}
Let $D_n$ denote the number of derangements (permutations without
fixed points) in $S_n$.  
Then
\begin{align*}
d_{\lceil n/2 \rceil} (n) &= D_{\lfloor n/2 \rfloor}, n \geq 4\\
d^*_{\lceil (n+1)/2 \rceil} (n) &= D_{\lfloor (n-1)/2 \rfloor}, n \geq 5. 
\end{align*}
\end{conj}

\begin{figure}[ht]
\centering
\includegraphics[height=3.5in]{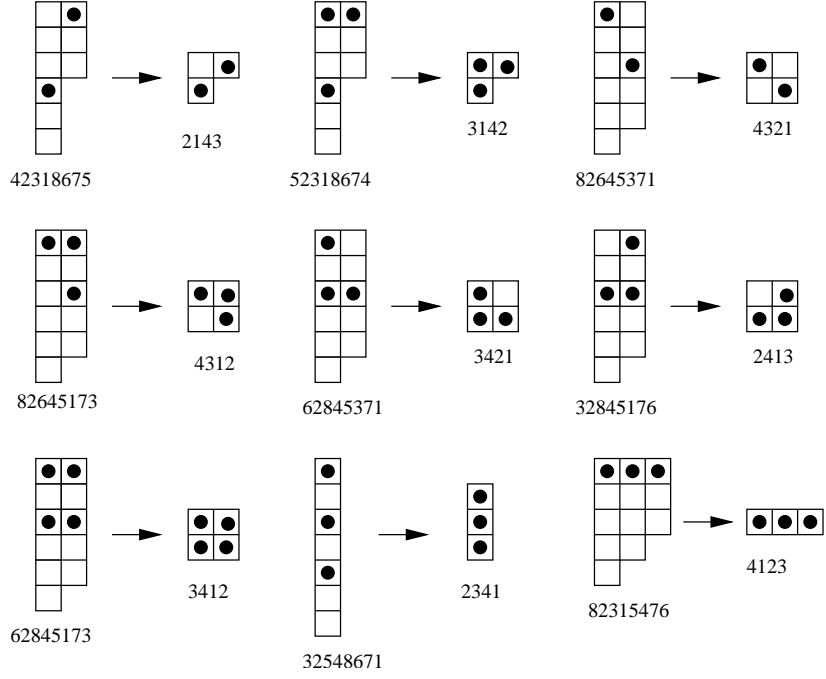}
\caption{The bijection $\Psi$ for $n=8$}
\label{Example}
\end{figure}

In this note we will give
two proofs of 
his conjecture relating the number of derangements to the number
of alternating permutations with 
the maximal number of fixed points.   Both proofs use the same bijection 
$\Psi$. 
The first proof works directly with permutations and shows that  
$\Psi$ is injective and surjective.  The second proof works with
permutation tableaux, certain tableaux
which are naturally in bijection with permutations, and explicitly 
constructs the inverse to $\Psi$.  The bijection (for alternating permutations)
is illustrated in 
Figure \ref{Example}, in terms of both permutations and 
permutation tableaux.

\section{The first proof}

As we will show subsequently, the main case that one needs to consider 
concerns alternating permutations on an even number of letters.

\begin{theorem}\label{Main}
For each nonnegative integer $m$
$$d_m(2m)=D_m.$$
\end{theorem}
\begin{proof}
Let $\ca_m$ denote the set of alternating
permutations of $[2m]$ with exactly $m$ fixed points
and let $\cd_m$ denote the set of derangements of $[m]$.
We shall define an explicit bijection $\Psi:\ca_m\to\cd_m$.

Let $\pi\in S_{2m}$ be an alternating permutation. We partition
the set $[2m]$ into $m$ two-element subsets $I_1,\ldots,I_m$
by setting $I_j=\{2j-1,2j\}$. As $\pi$ is alternating, $\pi(2j-1)>\pi(2j)$
and so $2j-1$ and $2j$ cannot both be fixed points of~$\pi$. Hence
$\pi$ has at most $m$ fixed points, and if $\pi\in\ca_m$ then
$\pi$ fixes exactly one point in each $I_j$ and moves the other point.
Write $I_j=\{a_j,b_j\}$ where $\pi(a_j)=a_j$ and $\pi(b_j)\ne b_j$.
Thus $\pi$ permutes the $b_j$ with no fixed points, that is there is
a unique derangement $\si\in\cd_m$ with $\pi(b_j)=b_{\si(j)}$.
Set $\Psi(\pi)=\si$.

To prove that $\Psi$ is a bijection we prove that it is injective and
surjective. Suppose $\pi\in\ca_m$ and $\Psi(\pi)=\si$. If $\si(j)>j$
then
$$\pi(b_j)=b_{\si(j)}\ge2\si(j)-1>2j\ge a_j=\pi(a_j).$$
As $\pi$ is decreasing on $I_j$
then $b_j=2j-1$ and $a_j=2j$. Similarly if $\si(j)<j$ then $a_j=2j-1$
and $b_j=2j$. Hence $\si$ determines the $a_j$ and $b_j$. Then as
$\pi(a_j)=a_j$ and $\pi(b_j)=b_{\si(j)}$ then $\si$ also determines~$\pi$.
Thus $\Psi$ is injective.

To prove $\Psi$ is surjective take $\si\in\cd_m$. Define
$$(a_j,b_j)=\left\{
\begin{array}{cl}
(2j-1,2j)&\textrm{if $\si(j)<j$,}\\
(2j,2j-1)&\textrm{if $\si(j)>j$.}
\end{array}
\right.$$
Then each element of $[m]$ is labelled as either an $a_j$ or a~$b_j$.
Define $\pi\in S_{2m}$ by $\pi(a_j)=a_j$ and $\pi(b_j)=b_{\si(j)}$.
As $\si$ is a derangement then $\pi$ has exactly $m$ fixed points.
We claim that $\pi$ is alternating. First of all if $\si(j)<j$ then
$$\pi(2j)=\pi(b_j)=b_{\si(j)}\le2\si(j)<2j-1=\pi(2j-1).$$
A similar argument works also when $\si(j)>j$.
We also need to prove that $\pi(2j)<\pi(2j+1)$ for $1\le j\le m-1$.
Now as $\pi(2j)<\pi(2j-1)$ and one of $\pi(2j)=2j$ and $\pi(2j-1)=2j-1$
holds then $\pi(2j)\le 2j$. Similarly $\pi(2j+1)\ge 2j+1$. Hence
$\pi(2j)<\pi(2j+1)$ and $\pi$ is alternating. It follows that
$\pi\in\ca_m$ and clearly $\Psi(\pi)=\si$. Hence $\Psi$
is a bijection and this concludes the proof.
\end{proof}

Regarding a permutation as the list of its values we can describe
the action of $\Psi$ in a simple way. Take a permutation $\pi\in\ca_m$
and first delete all fixed points. Then replace each number $k$ in the
remaining list by $\lceil{k/2}\rceil$. For example take
$\pi=52318674$. Deleting its fixed points gives $5184$, and
halving and rounding up each entry gives $\Psi(\pi)=3142$.

The remaining cases are simple consequences of this.

\begin{theorem}
For each positive integer $m$
$$d^*_{m+1}(2m)=D_{m-1}$$
and
$$d_m(2m-1)=d^*_m(2m-1)=D_{m-1}.$$
\end{theorem}
\begin{proof}
A reverse alternating permutation of $[2m]$ having $m+1$ fixed
points must fix $1$ and~$2m$. Also it restricts to an alternating
permutation of $\{2,\ldots,2m-1\}$ with $m-1$ fixed points and so
$d^*_{m+1}(2m)=d_{m-1}(2m-2)=D_{m-1}$.
Conjugating with the reversal permutation $\rho:j\mapsto 2m-j$
of $[2m-1]$ interchanges alternating permutations in $S_{2m-1}$ with $m$ fixed
points with reverse alternating permutations in $S_{2m-1}$ with
$m$ fixed points. Hence $d_m(2m-1)=d^*_m(2m-1)$. Also an alternating
permutation of $[2m-1]$ with $m$ fixed points must fix $2m-1$ and so
restricts to an alternating permutation of $[2m-2]$ with $m-1$
fixed points. Hence $d_m(2m-1)=d_{m-1}(2m-2)=D_{m-1}$.
\end{proof}

\section{Permutation Tableaux}

A {\em partition} $\lambda = (\lambda_1, \dots,
\lambda_k)$ is a weakly decreasing sequence of nonnegative
integers. For a partition $\lambda$, the
{\em Young diagram} $Y_\lambda$ of shape $\lambda$ is a left-justified
diagram of boxes, with $\lambda_i$ boxes in the $i$-th row.

A {\em permutation tableau} $\T$ \cite{Postnikov, SW} is a partition $\lambda$
together with a filling of each box of $Y_\lambda$ with either
a (black) dot or nothing such that the following holds:
\begin{enumerate}
\item Each column of $Y_\lambda$ contains at least one dot.
\item There is no empty box which has a dot above it in the same column
{\em and} a dot to its left in the same row.
\end{enumerate}

Figure \ref{fig-diagram} gives an example of a permutation tableau.

\psset{unit=.7pt, arrowsize=7pt, linewidth=1pt}
\psset{linecolor=blue}
\newgray{grayish}{.90}
\newrgbcolor{embgreen}{0 .5 0}
\def\vblack(#1, #2)#3{\cnode*[linecolor=black](#1, #2){3}{#3}}
\def\vwhite(#1,#2)#3{\cnode[linecolor=black,fillcolor=white,fillstyle=solid](#1,#2){3}{#3}}
\countdef\x=23
\countdef\y=24
\countdef\z=25
\countdef\t=26

\def\tbox(#1,#2)#3{
\x=#1 \y=#2
\multiply\x by 36
\multiply\y by 36
\z=\x \t=\y
\advance\z by 36
\advance\t by 36
\psline(\x,\y)(\x,\t)(\z,\t)(\z,\y)(\x,\y)
\advance\x by 18
\advance\y by 18
\rput(\x,\y){{\bf #3}}}

\def\evbox(#1,#2)#3{
\x=#1 \y=#2
\multiply\x by 36
\multiply\y by 36
\advance\x by 10
\advance\y by 18
\rput(\x,\y){{\bf #3}}}

\def\ehbox(#1,#2)#3{
\x=#1 \y=#2
\multiply\x by 36
\multiply\y by 36
\advance\x by 18
\advance\y by -12
\rput(\x,\y){{\bf #3}}}

\def\contournumbers{
\evbox(3,0){5}
\evbox(4,1){3}
\evbox(4,2){2}
\evbox(4,3){1}
\ehbox(0,0){8}
\ehbox(1,0){7}
\ehbox(2,0){6}
\ehbox(3,1){4}
}

\def\bdot{\pscircle*{1.5mm}}

\begin{figure}
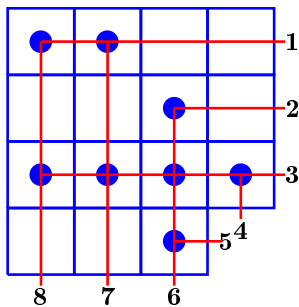

\pspicture(10,-10)(144,160)
\tbox(0,0){}
\tbox(1,0){}
\tbox(2,0){\bdot}

\tbox(0,1){\bdot}
\tbox(1,1){\bdot}
\tbox(2,1){\bdot}
\tbox(3,1){\bdot}

\tbox(0,2){}
\tbox(1,2){}
\tbox(2,2){\bdot}
\tbox(3,2){}

\tbox(0,3){\bdot}
\tbox(1,3){\bdot}
\tbox(2,3){}
\tbox(3,3){}

\contournumbers

\psset{linecolor=red}

\psline(150,126)(18,126)
\psline(18,126)(18,-6)

\psline(150,54)(18,54)
\psline(54,126)(54,-6)

\psline(90,90)(90,-6)
\psline(90,18)(116,18)

\psline(90,90)(150,90)
\psline(126,54)(126,30)
\endpspicture
\caption{\label{fig-diagram} The diagram of a tableau.}  
\end{figure}


We now recall the bijection $\Phi$ from permutation
tableaux to permutations.  More precisely, $\Phi$ is a bijection from
the set of permutation tableaux with $k$ rows and $n-k$ columns
 to permutations in the symmetric group
$\cls_n$ with $k$ weak excedances.  Here, a {\em weak excedance} of a
permutation $\pi$ is a value $\pi(i)$ such that $\pi(i) \geq i$.  
If we 
define the {\it semiperimeter} of a tableau to be the number of 
columns plus the 
number of rows, 
then $\Phi$ maps the set of permutation tableaux of semiperimeter $n$
to permutations in $\cls_n$.

We define the {\em diagram}
$D(\T)$ associated with $\T$ as follows.  
Label the edges of the northeast border of the Young diagram with 
the numbers from $1$ to $n$, as in Figure \ref{fig-diagram}.
From each black dot $v$, draw an
edge to the east and an edge to the south; each such edge should
connect $v$ to either a closest vertex in the same row or column, or
to one of the labels from $1$ to $n$.  The resulting picture is the
{\em diagram} $D(\T)$.  

We now define the permutation $\pi = \Phi(\T)$ via the following
procedure.  For each $i \in \{1, \dots , n\}$, find the corresponding
position on $D(\T)$ which is labeled by~$i$.  If the label $i$ is on a
vertical step of $P$, start from this position and travel straight
west as far as possible on edges of $D(\T)$. Then, take a ``zig-zag''
path southeast, by traveling on edges of $D(\T)$ south and east and
turning at each opportunity (i.e. at each new vertex).  This path will
terminate at some label $j\ge i $, and we let $\pi(i) = j$.  If $i$ is
not connected to any edge
then we set $\pi(i)=i$.  Similarly, if the label $i$
is on a horizontal step of~$P$, start from this position and travel
north as far as possible on edges of $D(\T)$. Then, as before, take a
zig-zag path south-east, by traveling on edges of $D(\T)$ east and
south, and turning at each opportunity.  This path will terminate at
some label $j<i$, and we let $\pi(i) = j$.

\begin{figure}
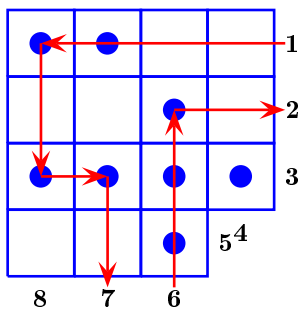

\pspicture(10,-10)(144,160)
\tbox(0,0){}
\tbox(1,0){}
\tbox(2,0){\bdot}

\tbox(0,1){\bdot}
\tbox(1,1){\bdot}
\tbox(2,1){\bdot}
\tbox(3,1){\bdot}

\tbox(0,2){}
\tbox(1,2){}
\tbox(2,2){\bdot}
\tbox(3,2){}

\tbox(0,3){\bdot}
\tbox(1,3){\bdot}
\tbox(2,3){}
\tbox(3,3){}

\contournumbers

\psset{linecolor=red}

\psline{->}(150,126)(18,126)
\psline{->}(18,126)(18,54)
\psline{->}(18,54)(54,54)
\psline{->}(54,54)(54,-6)

\psline{->}(90,-6)(90,90)
\psline{->}(90,90)(150,90)

\endpspicture
\caption{\label{zigzag}The paths taken by 1 and 6: $\pi(1)=7$, $\pi(6)=2$.}
\end{figure}

See Figure \ref{zigzag} for a picture of the path taken by $i$.

\begin{example}
If $\T$ is the permutation tableau whose diagram is
given in Figures \ref{fig-diagram} and \ref{zigzag}, then $\Phi(\T) =
74836215$.
\end{example}

The following lemma is clear from the construction above.
\begin{lemma}\label{fp}\cite{SW}
The positions of the weak excedances of $\pi = \Phi (\T)$ are precisely the
labels on the vertical edges of $P$.  The positions of non-excedances
of $\pi$ are
precisely the labels on the horizontal edges of $P$.  
Furthermore, in $\Phi(\T)$, the letter $i$ is a fixed point if and only if 
the row in $\T$ whose right hand edge is labeled by $i$ does not contain
any dots (is ``empty").
\end{lemma}


\section{The second proof}

Let $AT_n(2n)$ denote the set of permutation tableaux corresponding
to alternating permutations in $S_{2n}$ with $n$ fixed points, and let
$DT(n)$ denote the set of permutation tableaux corresponding to 
derangements in $S_n$.
We will describe $\Psi$ on the level of tableaux
and prove that it is a bijection between $AT_n(2n)$ and $DT(n)$.

By Lemma \ref{fp}, the set $DT(n)$ consists of 
permutation tableaux with semiperimeter $n$ such that 
no row is empty (i.e.\ 
each row contains at least one black dot).
We define $\Psi$ on the set $AT_n(2n)$: it acts on 
a tableau $\T$ by deleting all empty rows.  See Figure 
\ref{Example}.  Clearly
$\Psi(\T) \in DT(n)$ for some $n$. 

\begin{figure}[ht]
\centering
\includegraphics[height=1.5in]{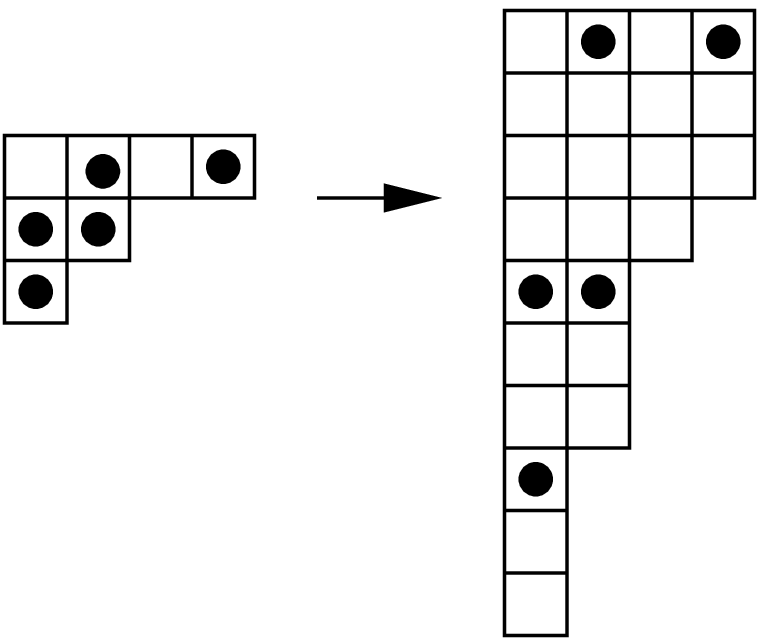}
\caption{}
\label{gExample}
\end{figure}

We now define a map $\Theta$ 
on $DT(n)$ (illustrated in Figure \ref{gExample}), which will be the inverse of $\Psi$, and acts
by inserting precisely $n$ empty rows into the tableau $\T$.
More specifically, if $\T \in DT(n)$ 
with partition shape $(\lambda_1, \lambda_2, \dots , \lambda_k)$, 
then $\Theta(\T)$ is the tableau which results 
after performing the following algorithm:
\begin{itemize}
\item For every $i$ such that $1 \leq i \leq k-1$, 
insert $\lambda_i - \lambda_{i+1}+1$ empty rows
 between the $i$th and $(i+1)$st rows of $\T$, whose rows lengths are
$(\lambda_i, \lambda_i, \lambda_i - 1, \lambda_i -2, \lambda_i - 3, \dots ,
   \lambda_{i+1}+1)$.
\item Insert $\lambda_k+1$ empty rows after the $k$th (last) row of $\T$, 
  whose row lengths are $(\lambda_k, \lambda_k, \lambda_k -1, \lambda_k - 2, \dots , 2,1)$.
\end{itemize}

The following lemma is obvious.

\begin{lemma}\label{obvious}
Let $\pi$ be an alternating permutation.
If $k$ consecutive entries 
are fixed points,
then $k \leq 2$. Furthermore, if two consecutive entries
$a_i$ and $a_{i+1}$ are fixed points, then $a_i < a_{i+1}$ and $i$
is even.
\end{lemma}

We now give a second proof of Theorem \ref{Main}.

\begin{proof}
We need to prove that $\Psi$ gives a bijection from
$AT_n(2n)$
to $DT(n)$, 
whose inverse is $\Theta$.  It is obvious that $\Theta$ 
is an injection, that 
$\Psi \circ \Theta$ is the identity, and that
if $\T$ is a tableau corresponding to a derangement
then $\Theta(\T)$ has $n$ fixed points.  
However, we need to show
that $\Theta(\T)$ is alternating and  that
$\Psi$ is an injection.

We will first show that $\Psi$ is an injection.
Suppose we are trying to insert $n$
empty rows into $\T$, so as to get a tableau $\T'$ corresponding
to an alternating permutation $\pi$.  Clearly the first row of the resulting
tableau $\T'$ cannot be empty, because then $\pi(1)=1$ violates
the requirement than $\pi(1)>\pi(2)$.
Additionally, by Lemma \ref{obvious}, 
we cannot have three consecutive 
empty rows of the same length.

Consider two consecutive rows in $\T$ which have lengths 
$t \geq s$.  Let us try to maximize the number
of empty rows that we can insert between these rows, such that 
the resulting permutation $\pi$ will be alternating.  Suppose these rows
have lengths $r_1, r_2, \dots , r_m$, with 
$t \geq r_1 \geq r_2 \geq \dots \geq r_m \geq s$.  
As before, there cannot be three consecutive
$r_i$'s which are equal.  Furthermore, note that after the rows of 
length $t$, there cannot be two empty rows of the same
length.  Indeed, if there were two such rows (say with vertical 
edge labels $j$ and $j+1$) then since $j-1$ is the label of a 
{\it horizontal} step of the tableau -- i.e.\ the position of 
a nonexcedance -- we would have 
$\pi(j-1) < j-1, \pi(j)=j$, and $\pi(j)=j+1$.  This violates the 
requirement that $\pi$ be alternating.  

So the best we could hope for is to have empty rows of lengths
$t, t, t - 1, t - 2, \dots , 
s + 1, s$.  But in fact we cannot have 
an empty row of length $s$.  If we did, then if 
$j$ is the label of the vertical step of this row, then we would
have $\pi(j-1) < j-1$ (since $j-1$ would be the label of a horizontal
step, hence a position of a nonexcedance).  Additionally, we would 
have $\pi(j) = j$ and $\pi(j+1) \geq j+1$, and these three entries
would violate the requirement that $\pi$ be alternating.

Therefore if we want to insert $n$ empty rows 
into $\T$ such that the resulting permutation is alternating, 
then our only choice is to use the algorithm 
described in the definition of $\Theta$.  Therefore
$\Psi$ is an injection.

Finally it remains to show that  
the permutation corresponding to $\Theta(\T)$ 
is alternating.
Consider two consecutive rows of lengths $\lambda_i$ and $\lambda_{i+1}$
in $\T$, and suppose that the vertical column labels corresponding
to these rows are $j$ and $k$.  We will show that the entries
of the corresponding permutation from positions $j$ to $k$ are alternating.
Since the row of length $\lambda_i$ has at least one dot, we have
$\pi(j) > j$.  Since $j+1$ and $j+2$ label the vertical steps 
of empty rows, we have $\pi(j+1)=j+1$ and $\pi(j+2)=j+2$.  
Since $j+3$ labels a horizontal step, $\pi(j+3)<j+3$ (and hence 
less than $j+2$ since $\pi(j+2)=j+2$).
Since $j+4$ labels a vertical step, $\pi(j+4)=j+4$.
Continuing in this fashion, we finally check that 
$\pi(k-2)=k-2$, $\pi(k-1)<k-1$ (hence less than $k-2$), and $\pi(k)>k$.
Clearly $\pi$ is alternating between the positions $j$ and $k$, 
and since the same argument can be applied to the remaining positions
of $\pi$, we are done.

\end{proof}

\end{document}